# MULTIPLE LOCAL WHITTLE ESTIMATION IN STATIONARY SYSTEMS

By P. M. Robinson[1]

*London School of Economics*

Moving from univariate to bivariate jointly dependent long-memory time series introduces a phase parameter $(\gamma)$, at the frequency of principal interest, zero; for short-memory series $\gamma = 0$ automatically. The latter case has also been stressed under long memory, along with the "fractional differencing" case $\gamma = (\delta_2 - \delta_1)\pi/2$, where $\delta_1$, $\delta_2$ are the memory parameters of the two series. We develop time domain conditions under which these are and are not relevant, and relate the consequent properties of cross-autocovariances to ones of the (possibly bilateral) moving average representation which, with martingale difference innovations of arbitrary dimension, is used in asymptotic theory for local Whittle parameter estimates depending on a single smoothing number. Incorporating also a regression parameter $(\beta)$ which, when nonzero, indicates cointegration, the consistency proof of these implicitly defined estimates is nonstandard due to the $\beta$ estimate converging faster than the others. We also establish joint asymptotic normality of the estimates, and indicate how this outcome can apply in statistical inference on several questions of interest. Issues of implemention are discussed, along with implications of knowing $\beta$ and of correct or incorrect specification of $\gamma$, and possible extensions to higher-dimensional systems and nonstationary series.

**1. Introduction.** In the analysis of long-memory time series, two major issues emerge in multivariate extension of univariate results. One is the possibility of cointegration, whereby one or more linear combinations of the (stationary or nonstationary) observables reduces memory. In general, rules of large sample inference based on a no-cointegration assumption are invalidated by cointegration, and vice versa. The literature on cointegration under long memory is dwarfed by that under autoregressive (AR) unit roots, but

Received June 2007; revised June 2007.

[1]Supported by ESRC Grant R000239936 and RES-062-23-0036.

has been developed in several directions recently. Another distinctive multivariate feature, which has attracted very little attention, is phase, essentially the argument in polar co-ordinate representation of the cross-spectrum. This is a particularly interesting issue in a "semiparametric" setting, where the spectral density matrix is modeled only near zero frequency. For a jointly covariance stationary short-memory process, this matrix is continuous at zero frequency; thus, since the quadrature spectrum (the imaginary part of the cross-spectrum) is an odd function, it, and thus the phase, are zero there. In long-memory series, on the other hand, where spectra diverge at zero frequency, the cross-spectrum is discontinuous there, and the phase need not be zero. In the literature, essentially two values for the phase have been considered, albeit rather implicitly, with little discussion of implications.

The present paper develops large sample statistical inference, in a possibly cointegrated system, with unknown phase. The formal results focus on a bivariate system, extension of our techniques for establishing asymptotic statistical theory to a system of arbitrary dimension being seemingly relatively straightforward, albeit introducing issues of specification and implementation, whose detailed treatment would be lengthy; we include a brief discussion. We also focus on covariance stationary observable series. This becomes a theoretical possibility when we switch from an AR unit root cointegration setting to a fractional one, and it has been of recent practical interest in financial time series analysis. We include, however, a brief discussion of possible nonstationary extensions.

Consider a bivariate jointly covariance stationary process $u_t = (u_{1t}, u_{2t})'$, having spectral density matrix $f_u(\lambda)$ that satisfies

$$(1.1) \qquad f_u(\lambda) \sim \Phi(\lambda; \alpha_0)^{-1} \Omega_0 \bar{\Phi}(\lambda; \alpha_0)^{-1} \qquad \text{as } \lambda \to 0,$$

$$(1.2) \qquad \Phi(\lambda; \alpha) = \text{diag}\{|\lambda|^{\delta_1}, |\lambda|^{\delta_2} e^{-i \, \text{sign}(\lambda)\gamma}\}, \qquad \lambda \in (-\pi, 0) \cup (0, \pi].$$

Here, $\alpha = (\gamma, \delta')'$ for $\delta = (\delta_1, \delta_2)'$, where $\gamma$, $\delta_1$ and $\delta_2$ are real-valued, $\gamma_0$ and $\delta_0 = (\delta_{01}, \delta_{02})'$ in $\alpha_0 = (\gamma_0, \delta_0')'$ are unknown, $\delta_{0i} \in [0, \frac{1}{2})$, $i = 1, 2$, $\Omega_0$ is an unknown $2 \times 2$ positive definite matrix, and the overbar indicates complex conjugation; the notation "$\sim$" in (1.1) means that for each element, the ratio of real/imaginary parts of the left and right sides tends to 1 (taking $0/0 = 1$).

From (1.1), $u_{it}$ is said to have memory (parameter) $\delta_{0i}$, its spectral density $f_i(\lambda)$ satisfying

$$f_i(\lambda) \sim \omega_{ii} |\lambda|^{-2\delta_{0i}} \qquad \text{as } \lambda \to 0, \ i = 1, 2,$$

where $\omega_{ij}$ is the $(i, j)$th element of $\Omega$. We deduce also that $u_{1t}, u_{2t}$ have cross-spectrum $f_{12}(\lambda)$ [the top right element of $f_u(\lambda)$] satisfying

$$(1.3) \qquad f_{12}(\lambda) \sim \omega_{12} |\lambda|^{-\chi_0} e^{-i \, \text{sign}(\lambda)\gamma_0} \qquad \text{as } \lambda \to 0,$$

where $\chi_0 = \delta_{02} + \delta_{01}$. Then (see, e.g., [4], page 302, [12], page 48) $\gamma_0$ is the phase between $u_{1t}, u_{2t}$ at $\lambda = 0$. There is no loss of generality in the restriction $\gamma_0 \in (-\pi, \pi]$. Thus the local approximation on the right of (1.3) is real-valued only if $\omega_{12} = 0$ and/or

$$\gamma_0 = 0. \tag{1.4}$$

To deduce another leading possibility, which applies to an extension of the fractional ARMA class, a general model for $f_u(\lambda)$ is

$$f_u(\lambda) = \Upsilon(\lambda; \alpha_0)^{-1} f_*(\lambda) \bar{\Upsilon}(\lambda; \alpha_0)^{-1}, \qquad \lambda \in (-\pi, 0) \cup (0, \pi], \tag{1.5}$$

where $\Upsilon(\lambda; \alpha) = \mathrm{diag}\{v(\lambda)^{\delta_1}, v(\lambda)^{\delta_2} e^{-i\,\mathrm{sign}(\lambda)\gamma}\}, v(\lambda) = (1 - e^{i\lambda}) e^{i\,\mathrm{sign}(\lambda)\pi/2}$ and $f_*(\lambda)$ is continuous and Hermitian positive definite at $\lambda = 0$. Since $v(\lambda) \sim |\lambda|$ as $\lambda \to 0$, (1.1) holds. On the other hand, with $\nu_0 = \delta_{02} - \delta_{01}$,

$$\gamma_0 = \frac{\pi}{2} \nu_0 \tag{1.6}$$

gives $\Upsilon(\lambda; \alpha_0) = \mathrm{diag}\{(1 - e^{i\lambda})^{\delta_{01}}, (1 - e^{i\lambda})^{\delta_{02}}\} e^{i\,\mathrm{sign}(\lambda)\delta_{01}\pi/2}$, so since the scalar factor has modulus 1, $u_t$ fractionally integrates an $I(0)$ process; if the latter is ARMA, $u_t$ is fractional ARMA. [Note that (1.6) reduces to (1.4) when $\delta_{01} = \delta_{02}$.] However, the fractional integration operator was originally motivated in a parametric framework [1], and in a semiparametric one there seems no overriding reason to fix $\gamma_0$. More generally, (1.1) with $\gamma_0 = (\delta_{02} - \delta_{01})c\pi/2$ can be shown to result from generalizing the fractional differencing filter $1 - e^{i\lambda}$ to $(1 - e^{i|\lambda|^{1/c}\,\mathrm{sign}(\lambda)})^c$, $c \neq 0$.

We can investigate the time domain implications of general $\gamma_0$. The proof of the following theorem is left to Section 5.

THEOREM 1. *Denoting* $r_{12}(j) = \mathrm{cov}(u_{1j}, u_{20})$, $j \in \mathbb{Z}$, *assume* $\chi_0 > 0$ *and, for* $(\kappa_+, \kappa_-) \neq (0, 0)$,

$$b_j = r_{12}(j) - \{\kappa_+ 1(j \geq 0) + \kappa_- 1(j < 0)\} |j|^{\chi_0 - 1} \tag{1.7}$$

*satisfies*

$$|b_j - b_{j+1}| \leq K |b_j|/(|j| + 1), \qquad b_j = o(|j|^{\chi_0 - 1}) \qquad as \ |j| \to \infty, \tag{1.8}$$

*where $K$ throughout denotes an arbitrarily large positive generic constant. Then* (1.3) *holds with*

$$\gamma_0 = \arctan\left\{\left(\frac{\kappa_+ - \kappa_-}{\kappa_+ + \kappa_-}\right) \tan \frac{\pi}{2} \chi_0\right\},$$
$$\omega_{12} = (\kappa_+ + \kappa_-)\Gamma(\chi_0) \cos(\pi\chi_0/2)/(2\pi \cos \gamma_0). \tag{1.9}$$

*In particular:*

$$\kappa_- = 0 \ is \ equivalent \ to \ \gamma_0 = \frac{\pi}{2}\chi_0, \qquad \omega_{12} = \kappa_+ \Gamma(\chi_0)/(2\pi), \tag{1.10}$$

$$\kappa_+ = 0 \ is \ equivalent \ to \ \gamma_0 = -\frac{\pi}{2}\chi_0, \qquad \omega_{12} = \kappa_- \Gamma(\chi_0)/(2\pi). \tag{1.11}$$

Solving (1.9) gives $\kappa_{\pm} = \pi\omega_{12}\sin(\pi\chi_0/2 \pm \gamma_0)/\Gamma(\chi_0)$. In view of (1.7) and the second part of (1.8), $r_{12}(j)$ dominates $r_{12}(-j)$ as $j \to \infty$ in (1.10), and vice versa in (1.11), while they decay at equal rates otherwise. The first part of (1.8) implies, with (1.7), an analogous condition for $r_{12}(j)$, which is satisfied by vector fractional ARMA processes. When $\kappa_+ = \kappa_-$ in (1.9), the power-law approximation is symmetric in $j$, and (1.4) results. On the other hand, (1.10) is a kind of weak causality ($u_2 \to u_1$) condition; it agrees with (1.6) only if $\delta_{01} = 0$. In general, the theorem indicates that any value of $\gamma_0$ is a possibility.

For the bivariate series $z_t = (y_t, x_t)'$, observed for $t = 1, \ldots, n$, consider the system

$$(1.12) \qquad B_0 z_t = u_t, \qquad t \in \mathbb{Z}, \qquad B_0 = \begin{pmatrix} 1 & -\beta_0 \\ 0 & 1 \end{pmatrix}$$

with $\beta_0$ unknown, so $u_{1t}$ is unobservable. When $\delta_{01} \geq \delta_{02}$, $\beta_0$ cannot be identified [from the spectral density matrix $f_z(\lambda)$ of $z_t$ near $\lambda = 0$] unless $\Omega_0$ is suitably restricted, for example, $\omega_{12}$ is known. When $\delta_{01} \neq \delta_{02}$, and $\beta_0 = 0$, $y_t$ and $x_t$ have unequal memories $\delta_{01}, \delta_{02}$, respectively. When $\delta_{01} < \delta_{02}$ and $\beta_0 \neq 0$, then both $x_t$ and $y_t$ have the same memory $\delta_{02}$, but the unobservable linear combination $u_{1t} = y_t - \beta_0 x_t$ has less memory, $\delta_{01}$, and $x_t$ and $y_t$ are said to be cointegrated. Both have a dominant common component with memory $\delta_{02}$, and so a dimensionality reduction is achievable:

$$(1.13) \qquad f_z(\lambda) \sim (\beta_0, 1)'(\beta_0, 1)\omega_{22}|\lambda|^{-2\delta_{02}}, \qquad \text{as } \lambda \to 0.$$

The right-hand side of (1.13) is singular, and the cointegrating error $u_{1t}$ has memory $\delta_{01}$. Included is the possibility that $\delta_{01} = 0$, when $u_{1t}$ has short memory. We focus on estimating $\theta_0 = (\beta_0, \alpha_0')'$ under

$$(1.14) \qquad 0 \leq \delta_{01} < \delta_{02} < \tfrac{1}{2},$$

covering cointegrated systems ($\beta_0 \neq 0$), and, for $\delta_{01} < \delta_{02}$, noncointegrated ones ($\beta_0 = 0$).

In [31] estimation of $\beta_0$ in (1.12) was discussed with $z_t$ exhibiting quite general forms of nonstationarity, and $u_{1t}$ being stationary or nonstationary. Reference [27] pointed out that cointegration is possible even when $z_t$ is stationary with long memory, as might be true of certain financial time series, say, and a number of references (e.g., [6, 23, 24]) have developed theory and applications in this setting. Financial time series are often very long, motivating reliance on only the "semiparametric," local, assumption (1.1). This justifies methods with only slow convergence rates, but a very large $n$ compensates. Faster rates are available in parametric models, for example when $u_t$ is a fractional ARMA process. However, if the ARMA component is misspecified, in that either the autoregressive (AR) or moving average (MA) orders are underspecified, or both are overspecified, all parameters will be inconsistently estimated. In [5] estimation of cointegrating subspaces in a semiparametric fractional context was studied. A recent parametric reference is [19].

We consider a narrow-band or local Whittle estimate $\hat{\theta} = (\hat{\beta}, \hat{\alpha}')' = (\hat{\beta}, \hat{\gamma}, \hat{\delta}_1, \hat{\delta}_2)'$ extending that for scalar long-memory series of [20], whose asymptotic properties were developed by [29], and further studied by and extended to nonstationary or noncointegrated multivariate systems by [18, 22, 26, 33, 34, 35]. References [36, 37] considered a version of it for cointegrated systems but with nonstationary fractional observables, while [24] has alternative results in the stationary case. We establish asymptotic properties of $\hat{\theta}$. For estimates that are only implicitly defined, a central limit theorem (CLT) is typically preceded by a consistency proof. This is more difficult to establish than usual because $\hat{\beta}$ converges faster than $\hat{\alpha}$. Consistency is usually established by showing that, after suitable normalization, the objective function converges uniformly in the parameter space to a limit which identifies all parameters and can thus be uniquely optimized. In multiparameter models this approach only works when all parameter estimates converge at the same rate. Additionally, as encountered by Robinson [29] in local Whittle estimation of the memory of a scalar series, our consistency result is insufficient to show that in the usual mean value theorem relations commencing the CLT proof, points on line segments between $\hat{\theta}$ and $\theta_0$ can be replaced to negligible effect by $\theta_0$; a slow convergence rate for $\hat{\delta}_1, \hat{\delta}_2$ is needed, and established using the stronger moment condition in any case required for the CLT.

The following section describes $\hat{\theta}$. Section 3 presents regularity conditions, a consistency result and CLT, and a small simulation study of finite-sample performance. Section 4 contains further discussion. Proofs are in Sections 5–8.

## 2. Local Whittle estimation.

For a generic vector $w_t$ define the periodogram matrix $I_w(\lambda) = n^{-1}(\sum_{t=1}^n w_t e^{it\lambda})(\sum_{t=1}^n w_t e^{-it\lambda})'$. Define the Fourier frequencies $\lambda_j = 2\pi j/n$, for integer $j$. In connection with (1.2) we allow some choice of "working model" for $f_u(\lambda)$ near $\lambda = 0$. Introduce

$$\Psi(\lambda; \alpha) = \text{diag}\{\psi(\lambda)^{\delta_1}, \psi(\lambda)^{\delta_2} e^{-i \, \text{sign}(\lambda)\gamma}\},$$

for a given complex-valued function $\psi(\lambda)$ such that $\psi(-\lambda) = \bar{\psi}(\lambda)$ and

$$(2.1) \qquad \psi(\lambda) - |\lambda| = o(1) \qquad \text{as } \lambda \to 0.$$

For example, $\psi(\lambda) = |\lambda|$ or $\upsilon(\lambda)$. Defining $A(\lambda; \theta) = \Psi(\lambda; \alpha) B I_z(\lambda) B' \bar{\Psi}(\lambda; \alpha)$, where $\theta = (\beta, \alpha')'$ and $B$ is defined as in (1.12) with $\beta_0$ replaced by $\beta$, consider the objective function

$$Q(\theta, \Omega) = \frac{1}{m} \sum_{j=1}^m [\log \det\{\Psi(\lambda_j; \alpha)^{-1} \Omega \bar{\Psi}(\lambda_j; \alpha)^{-1}\} + \text{tr}\{A(\lambda_j; \theta)\Omega^{-1}\}],$$

for $\Omega \in S$, the set of real positive definite $2 \times 2$ matrices, and an integer $m \in [1, n/2]$ which satisfies at least

$$(2.2) \qquad \frac{1}{m} + \frac{m}{n} \to 0 \qquad \text{as } n \to \infty.$$

The real function $Q$ is minimized over $S$ by $\hat{\Omega}(\theta) = \text{Re}\{m^{-1} \sum_{j=1}^{m} A(\lambda_j; \theta)\}$, leading to

$$R(\theta) = Q(\theta, \hat{\Omega}(\theta)) = \log \det\{\hat{\Omega}(\theta)\} - 2(\delta_1 + \delta_2) \frac{1}{m} \sum_{j=1}^{m} \log |\psi(\lambda_j)|.$$

Thus estimate $\theta_0$ by $\hat{\theta} = \arg\min_\Theta R(\theta)$, for a compact set $\Theta \in \mathbb{R}^4$ such that $\Theta = \Theta_\beta \times \Theta_\gamma \times \Theta_\delta$, with $\Theta_\beta$, $\Theta_\gamma$, $\Theta_\delta$ chosen as follows. Take $\Theta_\delta = \{\delta : -\eta_1 \le \delta_1 \le \delta_2 - \eta_2 \le \frac{1}{2} - \eta_2 - \eta_3\}$, where the $\eta_i$ are arbitrarily small positive numbers satisfying $0 < \eta_1 < \min(\eta_2, \eta_3)$, $\eta_2 + \eta_3 < \frac{1}{2}$; our consistency proof necessitates including a constraint corresponding to (1.14). We allow some $\delta_1 < 0$ because the CLT requires $\theta_0$ to be interior to $\Theta$, and we cover short memory, $\delta_{01} = 0$. We choose $\Theta_\gamma = [\eta_4 - \pi/2, \pi/2 - \eta_4]$ for $\eta_4 \in (0, \eta_3 - \eta_1)$, so $\gamma_0 \in \Theta_\gamma$ under (1.4) and (1.6). We can take $\Theta_\beta$ to be an arbitrarily large interval, possibly including $\{0\}$.

## 3. Asymptotic and finite-sample properties.

Existence of $f_u(\lambda)$ implies that for $p \ge 2$ we can find a $2 \times p$ matrix-valued function $C(\lambda)$ such that $C(-\lambda) = \overline{C}(\lambda)$ and

$$(3.1) \qquad f_u(\lambda) = C(\lambda) \bar{C}(\lambda)', \qquad \lambda \in (-\pi, \pi].$$

The representation (3.1) is familiar in case $p = 2$, but it is then obviously available for $p \ge 2$. Even when $p = 2$, $C(\lambda)$ is defined only up to postmultiplication by a unitary matrix, and when $p > 2$ the ambiguity is greater. From [12], page 61, existence of $f_u(\lambda)$ is equivalent to $u_t$ having representation

$$(3.2) \qquad u_t = Eu_t + \sum_{j \in \mathbb{Z}} C_j \varepsilon_{t-j}, \qquad t \in \mathbb{Z}, \ \sum_{j \in \mathbb{Z}} \|C_j\|^2 < \infty,$$

where $\{\varepsilon_t\}$ is a $p \times 1$ vector process such that $E\varepsilon_t = 0$, $E\varepsilon_t \varepsilon_t' = I_p$ (the $p \times p$ identity matrix), $E\varepsilon_s \varepsilon_t' = 0$, $s \ne t$, $s, t \in \mathbb{Z}$, $C_j = (2\pi)^{-1} \int_{-\pi}^{\pi} C(\lambda) e^{-ij\lambda} \, d\lambda$, and $\| \cdot \|$ is Euclidean norm. We will have to strengthen the conditions on $\varepsilon_t$ for asymptotic theory, but first discuss two other features of (3.2).

Moving average (MA) representations of long-memory time series models have typically been one-sided in particular $C_j = 0$, all $j < 0$, in (3.2), implying $u_t$ is purely nondeterministic (see, e.g., [11]). (An exception is [8], which considers a parametric model.) With Assumption A2, and the stronger Assumption B2 below for central limit theory, a one-sided representation was assumed in [29] in asymptotic theory for local Whittle estimation of memory parameter estimation, and subsequently by a number of authors in extensions of this work. On the other hand, since the basic quantity modeled is the spectral density matrix, rather than the process itself, there is no essential reason to impose one-sidedness. Indeed, going back to the earlier literature one can find repeated examples of bilateral representations in time series asymptotics (e.g., [2, 12, 25]). More recently, such representations have been employed to model specific (non-Gaussian, short-memory) phenomena (see, e.g., [3, 21], as well as examples in the electrical engineering literature, say). Our main motivation for allowing a bilateral representation here is to indicate its ability to yield any phase under long memory.

THEOREM 2. *Let (3.2) hold with $\{\varepsilon_t\}$ satisfying the conditions that follow it, and, denoting the $(k, \ell)$th element of $C_j$ by $c_{jk\ell}$, let*

$$g_{jk\ell} = c_{jk\ell} - \{\xi_{+k\ell}1(j \geq 0) + \xi_{-k\ell}1(j < 0)\}|j|^{\delta_{0k}-1}$$

*satisfy*

$$|g_{jk\ell} - g_{j+1,k\ell}| \leq K|g_{jk\ell}|/(|j|+1), \qquad g_{jk\ell} = o(|j|^{\delta_{0k}-1}) \qquad as \ |j| \to \infty,$$

*for constants $\xi_{+k\ell}, \xi_{-k\ell}$, $k = 1, 2$ and $\ell = 1, \ldots, p$. Then (1.7) and (1.8) of Theorem 1 hold with*

$$\kappa_+ = \xi'_{+1}\xi_{+2}B(1-\chi_0, \delta_{02}) + \xi'_{+1}\xi_{-2}B(\delta_{01}, \delta_{02}) + \xi'_{-1}\xi_{-2}B(1-\chi_0, \delta_{01}),$$

$$\kappa_- = \xi'_{+1}\xi_{+2}B(1-\chi_0, \delta_{01}) + \xi'_{+1}\xi_{-2}B(\delta_{02}, \delta_{01}) + \xi'_{-1}\xi_{-2}B(1-\chi_0, \delta_{02}),$$

*where $\xi_{+\ell}, \xi_{-\ell}$ are $p \times 1$ vectors with $k$th elements $\xi_{+k\ell}, \xi_{-k\ell}$, respectively.*

Section 6 contains a proof sketch. When $\xi_{-1} = \xi_{-2} = 0$, so that $u_t$ is purely nondeterministic, the relation $\Gamma(x)\Gamma(1-x) = -\pi \csc(\pi z)$ and trigonometric addition formulae may be shown to give (1.6), to extend the known results for fractional ARMA models. On the other hand, [6, 22, 23, 24] consider purely nondeterministic long-memory vector sequences with zero phases, (1.4), and we do not know of $C_j$ satisfying this prescription. However, the power-law decay of MA coefficients is only a sufficient condition for power-law spectral behavior. When $\xi_{+1} = \xi_{+2} = 0$, so $u_t$ has a one-sided forward representation, then $\gamma_0 = -\nu_0\pi/2$, the negative of (1.6), and the theorem indicates that for bilateral models $\gamma_0$ can take any value, which depends on the $\xi_{+\ell}, \xi_{-\ell}$ as well as the $\delta_{0i}$.

Another difference from the earlier references where MA representations are used in asymptotic theory for local Whittle estimates is in the allowance for rectangular, not necessarily square, $C_j$ in (3.2), and thus $u_t$ generated by shocks of higher dimension than the bivariate observable. Note that the equivalence property mentioned when introducing (3.2) is lost when $\varepsilon_t$ satisfies stronger assumptions, as in Assumption A2 below, but some generality can be recouped by the allowance for $p > 2$. This is natural if $x_t$, $y_t$ are seen as just two of a vector of related observations that are analyzed pairwise. It is also natural if (1.12) is viewed as a consequence of component models for $x_t$, $y_t$, namely $x_t = a_t + b_t$, $y_t = \beta_0 a_t + c_t$, where $a_t, b_t, c_t$ are unobservable sequences such that $a_t$ has memory $\delta_{02}$ and $u_{1t} = c_t - \beta_0 b_t$ has memory $\delta_{01}$; if the memories of $b_t$ and $c_t$ differ, then $b$ in Assumptions B1, B3 and B5 below is restricted. We can allow $(a_t, b_t, c_t)$ to have a nonsingular spectral density matrix by choosing $p \geq 3$ in (3.1). Note that $x_t$ and $y_t$ might themselves be instantaneous nonlinear functions of raw series $X_t$, $Y_t$, where $Y_t$ and $X_t$ are nonlinearly related, for example (in view of evidence of stationary long memory and cointegration in nonlinear functions of financial time series, see, e.g., [6]), logged squares, with $X_t$, $Y_t$ generated by long-memory stochastic volatility models, $X_t = A_t B_t$, $Y_t = A_t^\beta C_t$, where $A_t = e^{a_t}$, $B_t = e^{b_t}$, $C_t = e^{c_t}$.

We introduce the following assumptions for our consistency result.

ASSUMPTION A1.    Property (1.1) holds, where $u_t$ is covariance stationary, and for $C(\lambda)$ in (3.1),

$$(3.3) \qquad \Phi(\lambda; \alpha_0) C(\lambda) - P = o(1) \qquad \text{as } \lambda \to 0+,$$

where the real $2 \times p$ matrix $P$ satisfies $PP' = \Omega_0$, and $C(\lambda)$ is differentiable in a neighborhood of $\lambda = 0$, satisfying there

$$(3.4) \qquad \Phi(\lambda; \alpha_0) \frac{d}{d\lambda} C(\lambda) = O(\lambda^{-1}) \qquad \text{as } \lambda \to 0+.$$

ASSUMPTION A2.    $\{\varepsilon_t\}$ in (3.2) satisfy also $E(\varepsilon_t | \mathcal{F}_{t-1}) = E(\varepsilon_t), E(\varepsilon_t \varepsilon_t' | \mathcal{F}_{t-1}) = E(\varepsilon_t \varepsilon_t')$, a.s., $t \in \mathbb{Z}$, where $\mathcal{F}_t$ is the $\sigma$-field of events generated by $\varepsilon_s$, $s \leq t$, and also $P(\varepsilon_t' \varepsilon_t > \eta) \leq K P(X > \eta)$ for all $\eta > 0$ for some scalar nonnegative random variable $X$ such that $EX < \infty$.

ASSUMPTION A3.    Property (2.1) holds.

ASSUMPTION A4.    $\theta_0 \in \Theta$.

ASSUMPTION A5.    Property (2.2) holds.

ASSUMPTION A6.

$$(3.5) \qquad 0 < |\omega_{12}| < (\omega_{11} \omega_{22})^{1/2}.$$

Assumption A6 on the one hand implies $\Omega_0$ is positive definite, and on the other rules out

$$(3.6) \qquad \omega_{12} = 0,$$

when $u_{1t}, u_{2t}$ are incoherent at $\lambda = 0$ [cf. (1.3)]. Under (3.6) $\gamma_0$ is unidentifiable. We subsequently discuss related problems in which $\gamma_0$ is known and (3.6) is permitted. It could be covered in our theorems with extra detail, but while (3.6) is milder than the time domain orthogonality condition $r_{12}(j) = 0$, $j \in \mathbb{Z}$, it is less usual in the cointegration setting than (3.5), which tends to treat observables as jointly dependent. Assumption A1 implies (1.1), and this and other conditions are natural extensions or modifications of ones in [22, 29, 33].

THEOREM 3. *Let Assumptions A1–A6 hold. Then*

$$\hat{\alpha} \to_p \alpha_0, \hat{\beta} = \beta_0 + o_p\left(\left(\frac{m}{n}\right)^{\nu_0}\right) \qquad as\ n \to \infty.$$

To prove asymptotic normality we introduce the following assumptions.

ASSUMPTION B1. Assumption A1 holds, with the right-hand side of (3.3) replaced by $O(\lambda^b)$, for some $b \in (0, 2]$.

ASSUMPTION B2. Assumption A2 holds, with also the elements of $\varepsilon_t$ having a.s. constant third and fourth moments and cross-moments, conditional on $\mathcal{F}_{t-1}$.

ASSUMPTION B3. Property (2.1) holds for all $\gamma \in \Theta_\gamma$, after replacing its right-hand side by $O(\lambda^b)$, $b \in (0, 2]$.

ASSUMPTION B4. $\theta_0$ is an interior point of $\Theta$.

ASSUMPTION B5. For any $C < \infty$

$$(3.7) \qquad \frac{(\log m)^2 m^{1+2b}}{n^{2b}} + \frac{(\log n)^C}{m} \to 0 \qquad as\ n \to \infty.$$

The extensions of the previous conditions are similar to ones in earlier literature, the requirement $(\log n)^C/m \to 0$ coping, as in [33], with the fact that $\log n$ terms are not eliminated at the outset when $\psi(\lambda) = |\lambda|$. Define by $\Sigma$ the symmetric $4 \times 4$ matrix with $(k, \ell)$th element $\sigma_{k\ell}$, given by $\sigma_{11} = 2\mu\{(1 - 2\nu_0)^{-1} - (1 - \nu_0)^{-2}\cos^2(\gamma_0)\}\omega_{22}/\omega_{11}$, $\sigma_{12} = -2\mu(1 - \nu_0)^{-1}\sin(\gamma_0)(\omega_{12}/\omega_{11})$, $\sigma_{13} = 2\mu\nu_0(1 - \nu_0)^{-2}\cos(\gamma_0)\omega_{12}/\omega_{11}$, $\sigma_{14} = -2\mu\nu_0(1 - \nu_0)^{-2}\cos(\gamma_0)\omega_{22}/\omega_{11}$, $\sigma_{22} = -\sigma_{34} = 2\mu\rho^2$, $\sigma_{23} = \sigma_{24} = 0$, $\sigma_{33} = \sigma_{44} = 4 + \sigma_{34}$, where $\mu = (1 - \rho^2)^{-1}$, $\rho = \omega_{12}/(\omega_{11}\omega_{22})^{1/2}$. Write $\Delta_n = \text{diag}\{\lambda_m^{-\nu_0}, 1, 1, 1\}$ and let $N_k$ denote a $k$-variate normal variate.

TABLE 1
*Frequency of Wald test rejections, nominal 5% level*

| | | | $n = 128$ | | | | $n = 512$ | | | | $n = 2048$ | | | |
|---|---|---|---|---|---|---|---|---|---|---|---|---|---|---|
| $\delta_{01}$ | $\delta_{02}$ | $\rho$ | $m$ | $\beta$ | $\gamma$ | $\delta_1$ | $m$ | $\beta$ | $\gamma$ | $\delta_1$ | $m$ | $\beta$ | $\gamma$ | $\delta_1$ |
| 0.05 | 0.45 | 0.75 | 13 | 95.0 | 8.6 | 18.3 | 32 | 99.6 | 8.4 | 25.0 | 81 | 100 | 6.1 | 38.9 |
| 0.05 | 0.45 | 0.75 | 25 | 93.5 | 6.0 | 59.0 | 64 | 99.9 | 5.5 | 76.5 | 161 | 100 | 3.5 | 83.9 |
| 0.05 | 0.45 | 0.75 | 51 | 69.8 | 6.0 | 99.5 | 128 | 96.3 | 4.5 | 100 | 323 | 100 | 6.4 | 100 |
| 0.05 | 0.45 | 0.9 | 13 | 97.3 | 5.5 | 18.1 | 32 | 99.8 | 6.3 | 32.9 | 81 | 100 | 5.0 | 52.9 |
| 0.05 | 0.45 | 0.9 | 25 | 96.4 | 4.1 | 61.5 | 64 | 99.9 | 3.5 | 82.7 | 161 | 100 | 4.2 | 93.2 |
| 0.05 | 0.45 | 0.9 | 51 | 84.1 | 2.3 | 98.9 | 128 | 99.9 | 4.4 | 100 | 323 | 100 | 11.0 | 100 |
| 0.2 | 0.3 | 0.75 | 13 | 92.5 | 16.8 | 40.6 | 32 | 94.4 | 21.7 | 66.8 | 81 | 95.6 | 15.2 | 94.9 |
| 0.2 | 0.3 | 0.75 | 25 | 89.7 | 12.0 | 88.6 | 64 | 92.6 | 17.0 | 99.1 | 161 | 98.0 | 12.9 | 100 |
| 0.2 | 0.3 | 0.75 | 51 | 90.6 | 4.9 | 100 | 128 | 93.3 | 7.9 | 100 | 323 | 99.6 | 11.0 | 100 |
| 0.2 | 0.3 | 0.9 | 13 | 91.9 | 15.7 | 41.7 | 32 | 93.3 | 16.1 | 73.1 | 81 | 98.7 | 12.0 | 97.3 |
| 0.2 | 0.3 | 0.9 | 25 | 88.8 | 10.5 | 91.5 | 64 | 95.8 | 12.8 | 99.8 | 161 | 99.8 | 9.0 | 100 |
| 0.2 | 0.3 | 0.9 | 51 | 91.0 | 5.5 | 100 | 128 | 98.0 | 7.8 | 100 | 323 | 100 | 6.1 | 100 |

THEOREM 4. *Let Assumptions B1–B5 and A6 hold. Then as $n \to \infty$*

$$m^{1/2}\Delta_n(\hat{\theta} - \theta_0) \xrightarrow{d} N_4(0, \Sigma^{-1}).$$

A consistent estimate $\hat{\Sigma}$ of $\Sigma$ is formed by plugging $\hat{\theta}$ in place of $\theta_0$, and elements of $\hat{\Omega}(\hat{\theta})$ for those of $\Omega_0$. After also replacing $\Delta_n$ by $\widehat{\Delta}_n = \mathrm{diag}\{\lambda_m^{\hat{\delta}_1 - \hat{\delta}_2}, 1, 1, 1\}$, we can form asymptotically valid confidence regions for $\theta_0$, and also test hypotheses of interest, such as the linear homogeneous restrictions $\beta_0 = 0$ "no-cointegration"; (1.4) "zero-phase"; (1.6) "purely non-deterministic"; $\gamma_0 = (\delta_{01} + \delta_{02})\pi/2$ "weak causality"; $\delta_{01} = 0$ "short-memory cointegrating error." A small Monte Carlo study of finite-sample performance was carried out along such lines. To satisfy (1.1), $u_t$ was generated from the fractional ARMA $\mathrm{diag}\{(1 - L)^{\delta_{01}}, (1 - L)^{\delta_{02}}\}(1 - 0.5L)u_t = R^{1/2}\varepsilon_t$, where $L$ is the lag operator, the $\varepsilon_t$ are bivariate normal, and $R$ has elements 1 and 4 down the main-diagonal and off-diagonal element $2\rho$. Thus $\gamma_0 = \nu_0\pi/2$ (and $\omega_{12} = 4\rho/\pi$). We took $\delta_0 = (0.05, 0.45)'$ and $(0.2, 0.3)'$, $\rho = 0.75$ and 0.9, $\beta_0 = 1$. On each of 1000 replications, $\hat{\theta}$ was computed for three values of $m$, $[n^{2/3}/2], [n^{2/3}], 2n^{2/3}$ in each of three sample sizes, $n = 128$, 512 and 2048. We employed $\psi(\lambda) = |\lambda|$ (so local misspecification was incurred), and $\eta_1 = 0.01$, $\eta_2 = \eta_3 = 0.02$, $\eta_4 = 0.005$, $\Theta_\beta = [-3, 3]$. Table 1 gives Wald test rejection frequencies, at nominal two-sided 5% level, for the hypotheses $\beta_0 = 0$ (under "$\beta$"), (1.6) (under "$\gamma$") and $\delta_{01} = 0$ (under "$\delta_1$").

The second hypothesis is true so that size is measured, while the others are false so that power is measured. When $\delta_0 = (0.2, 0.3)'$ the gap $\nu_0$ is very small (and hard to detect); here the test on $\gamma_0$ is clearly oversized, even for large $n$, though matters improve for large $m$, and for $\delta_0 = (0.05, \, 0.45)'$ the sizes are better on average, albeit variable. For the test on $\delta_{01}$ power is poor for the smallest $m$, especially but unsurprisingly when $\delta_{01} = 0.05$, but increases satisfactorily with both. Power for testing $\beta_0$ is mostly very high. Overall, it seems hard to draw firm conclusions about the effect of $\rho$, while a relatively large $m$ appears to work best. Our technical results can be readily adapted to justify score and pseudo-likelihood-ratio-type tests.

## 4. Discussion.

REMARK 1. Lack of block-diagonality in $\Sigma$ suggests that correctly fixing $\alpha$ in $R(\theta)$ or employing an estimate $\tilde{\alpha}$ which converges faster than $m^{1/2}$ gives an estimate, $\hat{\beta}(\alpha)$, say, that is more efficient than $\hat{\beta}$, satisfying $m^{1/2} \lambda_m^{-\nu_0} \{ \hat{\beta}(\alpha_0) - \beta_0 \} \overset{d}{\to} N_1(0, \sigma_{11}^{-1})$. Going even further, but assuming (1.6), [15] provided an even more precise estimate of $\beta_0$, having the same efficiency as one minimizing $Q(\theta, \Omega)$ after replacing $\alpha$ and $\Omega$ by known $\alpha_0$ and $\Omega_0$; this estimate has also the advantage of a closed form representation. However, the need to select more than one bandwidth number, and in other respects suitably design the estimate of $\alpha_0$, and possibly $\Omega_0$, presents some disadvantage.

REMARK 2. On the other hand, computationally simpler but less efficient estimates than $\hat{\beta}$ are available. Reference [27] suggested the narrow-band least squares estimate

$$(4.1) \qquad \tilde{\beta} = \text{Re} \left\{ \sum_{j=1}^{m} I_{yx}(\lambda_j) \right\} \Big/ \sum_{j=1}^{m} I_x(\lambda_j),$$

where $(I_{yx}(\lambda), I_x(\lambda))'$ makes up the second column of $I_z(\lambda)$, and showed it to be consistent under very similar conditions to some of those for Theorem 1; [32] showed it is $(n/m)^{\nu_0}$-consistent (cf. Theorem 1). It advantageously avoids estimating $\alpha_0$. Reference [6] showed $\tilde{\beta}$ to be $(n/m)^{\nu_0} m^{1/2}$-consistent and asymptotically normal under (3.6) and $\chi_0 < 1/2$; [23] gave analogous results for a weighted version of (4.1). Even when a CLT for $\tilde{\beta}$, or another simple estimate, is available, the limiting variance depends on $\alpha_0$. Under (3.5), [32] showed that $(n/m)^{\nu_0}(\tilde{\beta} - \beta_0)$ converges in probability to a nonzero constant, so no useful inferential result is available. Our $\hat{\beta}$ corrects the bias.

REMARK 3. Simpler estimates of other parameters are available. We can estimate $\delta_{01}$ and $\delta_{02}$ using univariate local Whittle (see, e.g., [20, 29]), bivariate log-periodogram [28] or bivariate local Whittle [22, 33] techniques, though such estimation of $\delta_{01}$ requires a preliminary estimate of $\beta_0$. Given a preliminary estimate $\tilde{\beta}$, a simple estimate of $\gamma_0$ is

$$\tilde{\gamma} = \arctan\left[\operatorname{Im}\left\{\sum_{j=1}^m s(\lambda_j)\right\} \Big/ \operatorname{Re}\left\{\sum_{j=1}^m s(\lambda_j)\right\}\right],$$

where $s(\lambda) = I_{yx}(\lambda) - \tilde{\beta} I_x(\lambda)$.

REMARK 4. When $\gamma_0 = 0$, $\Sigma$ is block-diagonal with respect to $\hat{\beta}, \hat{\delta}$ on the one hand and $\hat{\gamma}$ on the other. Treating $\gamma_0$ as an unknown parameter seems unique in a long-memory setting, and it is worth noting the effects of its prior misspecification. Suppose we fix $\gamma = \gamma^*$ in $R(\theta)$, and then minimize with respect to $\beta, \delta$. Denoting $\theta_0^* = (\beta_0, \gamma^*, \delta_{01}, \delta_{02})'$, arguments like those in the proofs of Theorems 1 and 2 give

$$\hat{\Omega}(\theta_0^*) \xrightarrow{p} \begin{bmatrix} \omega_{11} & \omega_{12}\cos(\gamma^* - \gamma_0) \\ \omega_{12}\cos(\gamma^* - \gamma_0) & \omega_{22} \end{bmatrix}.$$

Likewise, taking $a \sim_p b$ to mean $a/b \to_p 1$ element-wise, calculations in the proof of Theorem 4 give

$$\frac{\partial \hat{\Omega}(\theta_0^*)}{\partial \beta} \underset{p}{\sim} \frac{2\lambda_m^{-\nu_0}}{1 - \nu_0} \begin{bmatrix} 2\omega_{12}\cos(\gamma_0) & \omega_{22}\cos(\gamma^*) \\ \omega_{22}\cos(\gamma^*) & 0 \end{bmatrix}.$$

Thus from (8.3) in the proof of Theorem 4 below,

$$\frac{\partial R(\theta_0^*)}{\partial \beta} \underset{p}{\sim} \frac{2\lambda_m^{-\nu_0}}{1 - \nu_0} \frac{\omega_{12}\omega_{22}\sin(\gamma^* - \gamma_0)\sin(\gamma^*)}{\omega_{11}\omega_{22} - \omega_{12}^2\cos^2(\gamma^* - \gamma_0)}.$$

It is readily seen that $(\partial/\partial \delta_k) R(\theta_0^*) \to_p 0$, but due to the nondiagonal limiting structure of $(\partial^2/\partial\theta\,\partial\theta') R(\theta_0^*)$, it appears that unless $\gamma^* = \gamma_0$, or $\gamma^* = 0$, not only is the $\beta_0$ estimate only $(n/m)^{\nu_0}$-consistent but the $\delta_{0i}$ estimates are inconsistent. When $\gamma_0 \neq \gamma^* = 0$, these estimates are asymptotically normal but their limiting variance matrix is complicated, and depends on $\gamma_0$. Our discussion suggests a more serious cost to incorrectly fixing $\gamma^* \neq 0$, for example, when $\gamma$ is replaced by $\nu\pi/2$ in $Q(\theta, \Omega)$, where $\nu = \delta_2 - \delta_1$; see (1.6). However, it can also be inferred that such bias problems are absent under (3.6). There are two cases of potential interest. In one, (3.6) is assumed a priori, in the other it is not; in both $\gamma_0$ is specified. In both cases the estimates of $\beta_0$, $\delta_{01}$, $\delta_{02}$, after correct centering and normalization as in Theorem 4, converge to independent zero-mean normal variates, whose variances can be deduced from the formulae in $\Sigma$ in the latter case (which the CLT of [24] addresses).

REMARK 5. On the other hand, if $\beta_0$ is known (e.g., to be zero, where there is no cointegration) we can infer from Theorems 3 and 4 that after correct centering and $m^{1/2}$ normalization, the estimates of $\gamma_0$ and $\delta_0$ are asymptotically independent, with limiting variances given in the inverse of the matrix consisting of the last three rows and columns of $\Sigma$. In fact the consistency proof is much simpler than that of Theorem 3, and the results hold for $|\delta_{0i}| < \frac{1}{2}$, $i = 1, 2$, with $\Theta_\delta$ chosen suitably.

REMARK 6. Also in the known noncointegrated case $\beta_0 = 0$, consistent estimation of $\gamma_0$, as well as of $\delta_0$, is relevant in inference based on the sample mean $\overline{z} = (z_1 + \cdots + z_n)/n$. Under our conditions it may be shown that as $n \to \infty$

$$\text{diag}\{n^{1/2-\delta_{01}}, n^{1/2-\delta_{02}}\}(\overline{z} - Ez_1)$$
$$\xrightarrow{d} N_2(0, (2\pi\omega_{ij}\cos((i-j)\gamma_0)/(\Gamma(\delta_{0i}+\delta_{0j}+2)\cos(\pi(\delta_{0i}+\delta_{0j})/2)))),$$

where the $(i,j)$th element of the $2 \times 2$ variance matrix is indicated. In [30] inference was developed in which the $\omega_{ij}$ and $\delta_0$ are replaced by consistent estimates [better than log-$n$-consistent in case of $\delta_0$, for which (3.7) suffices] but assuming $\gamma_0$ satisfies (1.6). If this assumption is incorrect, a corresponding confidence ellipse would be inconsistent. This kind of issue does not arise under short memory $\delta_{01} = \delta_{02} = 0$, where the variance matrix is $2\pi f_u(0)$, and phase is bound to be zero.

REMARK 7. An earlier version of this paper employed a different phase parameterization, $\phi\nu$, in place of $\gamma$. This naturally covers (1.4) ($\phi_0 = 0$) and (1.6) ($\phi_0 = \pi/2$), but is less natural in general, in view of Theorems 1 and 2. It affects the form of $\Sigma$, in particular giving nonzero $\sigma_{23}$ and $\sigma_{24}$. As a consequence, when $\beta_0$ is known the limiting variance matrix for estimation of $\alpha_0$ is no longer block-diagonal (cf. Remark 5), while if $\phi$ is incorrectly specified to a nonzero value (e.g., $\pi/2$), $\delta_0$ is estimated inconsistently; in Remark 4, with $\gamma$ likewise misspecified, this was due to estimating $\beta_0$. On the other hand, with the $\phi\nu$ parameterization, [33] compared the cases when $\phi$ is correctly fixed at $\pi/2$, and when $\phi$ is correctly fixed at 0 (where the limit distribution is the same as in Remark 5), finding greater precision in the former.

REMARK 8. To construct approximate Newton iterations, given an $i$th iterate $\hat{\theta}^{(i)}$, $i \geq 1$, we can form $\hat{\Sigma}^{(i)}$ by plugging in $\hat{\theta}^{(i)}$ for $\theta_0$ in $\Sigma$, replacing elements of $\Omega$ by those of $\hat{\Omega}(\hat{\theta}^{(1)})$, and then compute $\hat{\theta}^{(i+1)} = \hat{\theta}^{(i)} - \hat{\Sigma}^{(i)-1}(\partial/\partial\theta)R(\hat{\theta}^{(i)})$. Choices for $\hat{\theta}^{(1)}$ include estimates described in Remarks 2 and 3. If $\hat{\theta}^{(1)}$ satisfies $m^{1/2}\Delta_n(\hat{\theta}^{(1)} - \theta_0) = O_p(1)$, then $\hat{\theta}^{(2)}$ has the properties of $\hat{\theta}$ in Theorem 4. If the initial $\beta_0$ estimate is only $(n/m)^{\nu_0}$-consistent, as is (4.1), $\hat{\theta}^{(i)}$ should satisfy Theorem 4 for some finite $i$ but determination of a minimal $i$ depends on hypothesizing a rate of increase for $m$ with $n$, and on the unknown $\nu_0$. If a smaller $m$ is used in (4.1) than in $R(\theta)$, assuming the former $m$ increases sufficiently slowly relative to the latter one can justify $i = 2$ even.

REMARK 9. With respect to choice of $m$, minimizing approximate mean squared error (MSE) of a given linear combination of $\hat{\theta}$ elements is complicated, especially as $\hat{\beta}$ converges faster than $\hat{\alpha}$. Though suboptimal, the minimum MSE rule (in scalar local Whittle estimation of memory) of [13] could be applied, most simply to the $x_t$ sequence (requiring preliminary estimation of $\delta_{02}$ and fixing $b$ in Assumption B1, say, to 2). As always, a minimum-MSE rate violates the assumption (here B5) that provides correct centering in the CLT, suggesting use of a smaller $m$. In univariate local Whittle memory estimation, with data tapering, Giraitis and Robinson [10] developed an $m$ that minimizes the error in the CLT, having rate $n^{b/(1+b)}$, which satisfies B5; with $b = 2$ this is the rate employed in the Monte Carlo. References [16] and [17] proposed data-dependent $m$ in univariate log-periodogram memory estimation. Full confidence cannot be placed in any automatic technique and it may be wise to employ a grid of $m$ values, and assess sensitivity; estimates for a given $m$ should be a good starting point for iterations with adjacent $m$.

REMARK 10. From Assumption B5, $\hat{\alpha}$, $\hat{\beta}$ converge slower than $n^{b/(1+2b)}$, $n^{1/2-(1/2-\nu_0)/(1+2b)}$, respectively, for example, $n^{2/5}$, $n^{(2+\nu_0)/5}$ for $b = 2$, while for all $b$ the rate of $\hat{\beta}$ approaches $n^{1/2}$ as $\nu_0 \to \frac{1}{2}$. This rate is best for estimates of all parameters if $f_u(\lambda)$, $\lambda \in (-\pi, \pi]$, is parametric (extending theory of [7, 9, 11, 14]). But misspecification of $f_u$ incurs inconsistent estimation of all parameters, and if $f_u$ involves additional parameters [over those in (1.1)], computational burden increases. The least squares estimate of $\beta_0$ is inconsistent when $u_{1t}$ and $u_{2t}$ are correlated (cf. Assumption A6).

REMARK 11. By analogy with the pseudospectrum of univariate nonstationary fractional series, we can define a pseudospectral density matrix [involving a phase parameter as in (1.1)] for vector series with one or more nonstationary elements. Integer differencing of both series will not change phase, and may produce the stationary setting of the present paper. Given uncertainty as to whether or not the data are nonstationary, or about the degree of nonstationarity, alternative methods, already employed to extend univariate local Whittle estimates (e.g., [34, 35]), should produce analogous asymptotic properties to those in Theorems 3 and 4, albeit perhaps with some variance inflation, so long as the gap between memory parameters is less than $\frac{1}{2}$ [as in (1.14)]. If this gap exceeds $\frac{1}{2}$ optimal estimates have a faster rate, and mixed normal asymptotics [15]. Reference [36] considered local Whittle estimation with a gap exceeding $\frac{1}{2}$, but the estimate of $\beta_0$ achieves a slower convergence rate than is attainable even by such simple estimates as (4.1) and least squares when also the sum of memory parameters exceeds 1.

REMARK 12. Another kind of extension concerns multivariate series $z_t$ of dimension $q > 2$. Reference [24] considers local Whittle estimation with $q \geq 2$ and a single cointegrating relation, though with phases correctly assumed to be zero, (3.6) assumed in the CLT, and a consistency proof which, like ours, takes $q = 2$. More generally, $q > 2$ raises the possibility that the number, $r < q$, of cointegrating relations exceeds 1. In (1.12), $B_0$ can be redefined by replacing the 1's in the diagonal by blocks $I_r$ and $I_{q-r}$, with $\beta_0$ now being an $r \times (q - r)$ matrix. Likewise in (1.1), (1.2) the dimension is extended to $q$, with, for $j \in [2, q]$, the $j$th diagonal element of $\Phi(\lambda; \alpha)$ now being $|\lambda|^{\delta_j} e^{-i \operatorname{sign}(\lambda)(\gamma_1 + \cdots + \gamma_{j-1})}$, with $\delta_i < \delta_j$ for $i \leq r$, $j > r$. Thus $\alpha = (\gamma_1, \ldots, \gamma_{q-1}, \delta_1, \ldots, \delta_q)'$ unless, to mitigate possible curse of dimensionality and additional computational challenge, prior restrictions are imposed, for example, $\delta_1 = \cdots = \delta_r$ and/or $\delta_{r+1} = \cdots = \delta_q$. Such constraints could imply some zero $\gamma_i$ even under fractional integration assumptions [cf. (1.6), which is zero for $\delta_{01} = \delta_{02}$], but in general they can be unrestricted. Prior restrictions on $\beta_0$ might also be imposed. Our methods can be straightforwardly extended to estimate the remaining, unknown, parameters. The techniques of proof of Theorems 3 and 4 also appear to extend, while Theorems 1 and 2 clearly remain relevant.

**5. Proof of Theorem 1.** From [38], page 186,

$$(5.1) \qquad \sum_{j=1}^{\infty} j^{\chi_0 - 1} e^{ij\lambda} = \Gamma(\chi_0) e^{i\pi\chi_0/2} \lambda^{-\chi_0} + O(1) \qquad \text{as } \lambda \to 0+.$$

For $\lambda \neq 0$, $\operatorname{mod}(2\pi)$,

$$f_{12}(\lambda) = (2\pi)^{-1} \left\{ r_{12}(0) + \kappa_+ \sum_{j=1}^{\infty} j^{\chi_0 - 1} e^{-ij\lambda} + \kappa_- \sum_{j=1}^{\infty} j^{\chi_0 - 1} e^{ij\lambda} + \sum_{|j|=1}^{\infty} b_j e^{-ij\lambda} \right\}.$$

The last term in braces is bounded by

$$\sum_{j=1}^{N} (|b_j| + |b_{-j}|) + \sum_{j=N+1}^{\infty} \{|b_j - b_{j+1}| + |b_j - b_{-j-1}|\} \left| \sum_{k=N}^{\infty} e^{-ik\lambda} \right|$$

$$\leq K\varepsilon (N^{\chi_0} + N^{\chi_0 - 1} |\lambda|^{-1}) = o(|\lambda|^{-\chi_0}) \qquad \text{as } \lambda \to 0,$$

where $\varepsilon > 0$ is arbitrary and we choose $N \sim |\lambda|^{-1}$. Thus from (5.1),

$$f_{12}(\lambda) \sim (2\pi)^{-1} (\kappa_+ e^{-i \operatorname{sign}(\lambda) \pi \chi_0/2} + \kappa_- e^{i \operatorname{sign}(\lambda) \pi \chi_0/2}) \Gamma(\chi_0) |\lambda|^{-\chi_0} \qquad \text{as } \lambda \to 0.$$

Then (1.9) is determined by inspection, and the remaining statements are straightforwardly verified.

**6. Proof of Theorem 2.** Take $j \geq 0$. With $c'_{ij}$ denoting the $i$th row of $C_j$, write

$$r_{12}(j) = \sum_{i=j+1}^{\infty} c'_{1i} c_{2,i-j} + \sum_{i=0}^{j} c'_{1i} c_{2,i-j} + \sum_{i=-\infty}^{-1} c'_{1i} c_{2,i-j}.$$

Each of the three terms on the right-hand side is dominated by contributions in which $c_{1i}, c_{2,i-j}$ are of order $|i|^{\delta_{01}-1}$ and $|i-j|^{\delta_{02}-1}$, respectively, the remainder terms involving products of these with the $g_{i1\ell}$, $g_{i-j,2\ell}$ and products of the latter. After integral approximation of the leading terms we write

$$
\begin{aligned}
(6.1) \quad r_{12}(j) = \Big\{ & \xi'_{+1}\xi_{+2} \int_0^{\infty} (1+x)^{\delta_{01}-1} x^{\delta_{02}-1}\, dx \\
& + \xi'_{+1}\xi_{-2} \int_0^1 x^{\delta_{01}-1}(1-x)^{\delta_{02}-1}\, dx \\
& + \xi'_{-1}\xi_{-2} \int_0^{\infty} x^{\delta_{01}-1}(1+x)^{\delta_{02}-1}\, dx \Big\} j^{\chi_0-1} + b_j.
\end{aligned}
$$

We omit the straightforward but lengthy proof that $b_j$ satisfies (1.7) and (1.8). It only remains to express the integrals in (6.1) as Beta functions. The method of proof for $j < 0$ is identical.

**7. Proof of Theorem 3.** We first give the proof with "$o$" replaced by "$O$" in the error bound for $\hat{\beta}$. For any $c > 0$ define neighborhoods $\mathcal{N}_\beta(c) = \{\beta : |\beta - \beta_0| < c\}$, $\mathcal{N}_\gamma(c) = \{\gamma : |\gamma - \gamma_0| < c\}$, $\mathcal{N}_\delta(c) = \{\delta : \|\delta - \delta_0\| < c\}$. Fix $\varepsilon > 0$ and define $\mathcal{N}(\varepsilon) = \mathcal{N}_\beta(\varepsilon^{-1}(m/n)^{\nu_0}) \times \mathcal{N}_\gamma(\varepsilon) \times \mathcal{N}_\delta(\varepsilon)$, $\overline{\mathcal{N}}(\varepsilon) = \Theta \setminus \mathcal{N}(\varepsilon)$. We have $P(\hat{\theta} \in \overline{\mathcal{N}}(\varepsilon)) \leq P(\inf_{\overline{\mathcal{N}}(\varepsilon)} \{R(\theta) - R(\theta_0)\})$. To show that this tends to zero we first decompose $R(\theta) - R(\theta_0)$. We omit the straightforward proof, using Assumption A3, that the effect of replacing $\Psi(\lambda; \alpha)$ by $\Phi(\lambda; \alpha)$, when they differ, is negligible, uniformly on $\overline{\mathcal{N}}(\varepsilon)$, and proceed as if $\Psi = \Phi$. Then

$$(7.1) \quad R(\theta) - R(\theta_0) = \log \det\{\hat{\Omega}(\theta)\hat{\Omega}(\theta_0)^{-1}\} - 2\sum_{i=1}^{2} \zeta_i \frac{1}{m}\sum_j \log \lambda_j,$$

where $\zeta_i = \delta_i - \delta_{0i}$ and $\sum_j$ means $\sum_{j=1}^{m}$. With $\Upsilon(\delta) = \mathrm{diag}\{(2\zeta_1 + 1)^{1/2}, (2\zeta_2 + 1)^{1/2}\}$, $\Xi(\theta) = \mathrm{diag}\{\lambda_m^{\zeta_1}, \lambda_m^{\zeta_2}\}$, $\hat{\Omega}^*(\theta) = \Xi(\delta)\hat{\Omega}(\theta)\Xi(\theta)$, write

$$R(\theta) - R(\theta_0) = \log \det\{\Upsilon(\delta)\hat{\Omega}^*(\theta)\Upsilon(\delta)\hat{\Omega}(\theta_0)^{-1}\} + u(\delta),$$

where $u(\delta) = \sum_{i=1}^{2}\{2\zeta_i - \log(2\zeta_i + 1) + 2\zeta_i(\log m - m^{-1}\sum_j \log j - 1)\}$. To decompose $\hat{\Omega}^*(\theta)$ denote $I_{uj} = I_u(\lambda_j)$, and deduce

$$BI_z(\lambda_j)B' = I_{uj} + (B_0 - B)I_{uj} + I_{uj}(B_0 - B)' + (B_0 - B)I_{uj}(B - B_0)'.$$

With the definitions $H_j = \Psi(\lambda_j; \alpha_0)I_{uj}\bar{\Psi}(\lambda_j; \alpha_0)$, $b_n(\beta) = \lambda_m^{-\nu_0}(\beta_0 - \beta)$, $\tau = \gamma - \gamma_0$, rearrangement gives

$$(7.2) \qquad \hat{\Omega}^*(\theta) = \hat{G}^{(1)}(\alpha) + b_n(\beta)\hat{G}^{(2)}(\alpha) + b_n^2(\beta)\hat{G}^{(3)}(\alpha),$$

where $\hat{G}^{(i)}(\alpha) = (\hat{g}_{k\ell}^{(i)})$, $\hat{g}_{kk}^{(1)} = m^{-1}\sum_j (j/m)^{2\zeta_k} h_{kkj}$ $(k = 1, 2)$, $\hat{g}_{12}^{(1)} = \hat{g}_{21}^{(1)} = (2m)^{-1}\sum_j (j/m)^{\zeta_1 + \zeta_2}(e^{i\tau}h_{12j} + e^{-i\tau}h_{21j})$, $\hat{g}_{11}^{(2)} = (2m)^{-1}\sum_j (j/m)^{\zeta_1 + \delta_1 - \delta_{02}} \times (e^{i\gamma_0}h_{21j} + e^{-i\gamma_0}h_{12j})$, $\hat{g}_{12}^{(2)} = \hat{g}_{21}^{(2)} = m^{-1}\sum_j (j/m)^{\delta_1 - \delta_{02} + \zeta_2}(\cos\gamma)h_{22j}$, $\hat{g}_{11}^{(3)} = m^{-1}\sum_j (j/m)^{2(\delta_1 - \delta_{02})}h_{22j}$, $\hat{g}_{22}^{(2)} = \hat{g}_{12}^{(3)} = \hat{g}_{21}^{(3)} = \hat{g}_{22}^{(3)} = 0$, suppressing reference to dependence on $\alpha$ in the $\hat{g}_{k\ell}^{(i)}$ and with $H_j = (h_{k\ell j})$. Defining

$$U_\alpha(\alpha) = \log\det\{\Upsilon(\delta)\hat{G}^{(1)}(\alpha)\Upsilon(\delta)\hat{G}^{(1)}(\alpha_0)^{-1}\} + u(\delta),$$

$$U_\beta(\theta) = \log\det\{\hat{\Omega}^*(\theta)\hat{G}^{(1)}(\alpha)^{-1}\},$$

we have $R(\theta) - R(\theta_0) = U_\alpha(\alpha) + U_\beta(\theta)$, since $\hat{\Omega}(\theta_0) = \hat{G}^{(1)}(\alpha_0)$. Writing $\bar{\mathcal{N}}_\beta(c) = \Theta_\beta \backslash \mathcal{N}_\beta(c)$, $\bar{\mathcal{N}}_\gamma(c) = \Theta_\gamma \backslash \mathcal{N}_\gamma(c)$, $\bar{\mathcal{N}}_\delta(c) = \Theta_\delta \backslash \mathcal{N}_\delta(c)$, and also $\Theta_\alpha = \Theta_\gamma \times \Theta_\delta$, $\bar{\mathcal{N}}_\alpha(c) = \{\bar{\mathcal{N}}_\gamma(c) \times \Theta_\delta\} \cup \{\Theta_\gamma \times \bar{\mathcal{N}}_\delta(c)\}$, it suffices to show that as $n \to \infty$

$$(7.3) \qquad P\left(\inf_{\bar{\mathcal{N}}_\alpha(\varepsilon)} U_\alpha(\alpha) \le 0\right) \to 0,$$

$$(7.4) \qquad P\left(\inf_{\bar{\mathcal{N}}_\beta(1/\varepsilon(n/m)^{\nu_0}) \times \Theta_\alpha} U_\beta(\theta) \le 0\right) \to 0.$$

Introduce the following population analogues of the $\hat{g}_{k\ell}^{(i)}$: $g_{kk}^{(1)} = \omega_{kk}(2\zeta_k + 1)^{-1}$ $(k = 1, 2)$, $g_{12}^{(1)} = g_{21}^{(1)} = (\zeta_1 + \zeta_2 + 1)^{-1}\omega_{12}\cos\tau$, $g_{11}^{(2)} = 2(\zeta_1 + \delta_1 - \delta_{02} + 1)^{-1}\omega_{12}\cos\gamma_0$, $g_{12}^{(2)} = g_{21}^{(2)} = (\delta_1 + \zeta_2 - \delta_{02} + 1)^{-1}\omega_{22}\cos\gamma$, $g_{11}^{(3)} = (2(\delta_1 - \delta_{02}) + 1)^{-1}\omega_{22}$, $g_{22}^{(2)} = g_{12}^{(3)} = g_{21}^{(3)} = g_{22}^{(3)} = 0$; write $G^{(i)}(\alpha) = (g_{k\ell}^{(i)})$.

To prove (7.3) observe that from the inequality $|\log(1 + x)| \le 2|x|$ for $|x| \le \frac{1}{2}$, and because $\bar{\mathcal{N}}_\alpha(\varepsilon) \subset \{\bar{\mathcal{N}}_\gamma(\varepsilon) \times \Theta_\delta\} \cup \bar{\mathcal{N}}_\delta(\varepsilon)$, it suffices (following a development like that in [22]) to show

$$(7.5) \qquad \sup_{\Theta_\alpha}\|\Upsilon(\delta)\{\hat{G}^{(1)}(\alpha) - G^{(1)}(\alpha)\}\Upsilon(\delta)\| \xrightarrow{p} 0,$$

$$(7.6) \qquad \sup_{\Theta_\alpha}\|\{\Upsilon(\delta)G^{(1)}(\alpha)\Upsilon(\delta)\}^{-1}\| < \infty,$$

$$(7.7) \qquad \inf_{\bar{\mathcal{N}}_\gamma(\varepsilon) \times \Theta_\delta} \log \det\{\Upsilon(\delta) G^{(1)}(\alpha) \Upsilon(\delta) G^{(1)}(\alpha_0)^{-1}\} > 0,$$

$$(7.8) \qquad \varliminf_{n \to \infty} \inf_{\bar{\mathcal{N}}_\delta(\varepsilon)} u(\delta) > 0.$$

We omit the details of (7.5) as these are now standard, mainly following the proof of Theorem 1 of [29], and multivariate extensions [22, 33]. Our model (3.2) is more general than those in such references in two respects, namely our allowance for a bilateral MA and for the dimension of $\varepsilon_t$ to exceed 2, but it is readily seen that neither extension materially affects the proof. The basic technique involves summation-by-parts (to deal with the uniformity) followed by approximation of the $H_j$ by the $P I_{\varepsilon j} P'$, where $I_{\varepsilon j} = I_\varepsilon(\lambda_j)$ (see [29]) and then approximating the consequent term in the $I_{\varepsilon j}$ by one in $\Omega_0$ (with only a second moment for $\varepsilon_t$ required for the latter step due to applying a law of large numbers for $L_1$ variables to the term in the $\varepsilon_t \varepsilon_t'$) and approximating sums of form $m^{-1} \sum_j (j/m)^a$ by $(1+a)^{-1}$ for $a > -1$. The most significant difference from earlier results is the presence of the general $\gamma, \gamma_0$, but this is easily handled in view of compactness of $\Theta_\gamma$. Likewise, (7.8) follows from the proof of Theorem 1 of [29], which used the inequalities

$$(7.9) \qquad \inf_{|x| > \varepsilon}\{x - \log(x+1)\} > \frac{\varepsilon^2}{6}, \qquad \left|\log m - m^{-1} \sum_j \log j - 1\right| \le K m^{-1}.$$

To prove (7.6) observe that

$$(7.10) \qquad \det\{\Upsilon(\delta) G^{(1)}(\alpha) \Upsilon(\delta)\} = \omega_{11} \omega_{22} - \omega_{12}^2 c(\delta) \cos^2 \tau,$$

where $c(\delta) = (2\zeta_1 + 1)(2\zeta_2 + 1)/(\zeta_1 + \zeta_2 + 1)^2$. It follows from the inequality $0 < 4xy \le (x+y)^2$, for $x, y > 0$, that $0 < c(\delta) \le 1$, and thus (7.10) $\ge \det(\Omega_0) > 0$.

To prove (7.7) note that

$$(7.11) \qquad \log \det\{\Upsilon(\delta) G^{(1)}(\alpha) \Upsilon(\delta) G^{(1)}(\alpha_0)^{-1}\} = \log\left\{\frac{1 - \rho^2 c(\delta) \cos^2 \tau}{1 - \rho^2}\right\}.$$

From $|\cos \tau| \le 1$, $|c(\delta)| \le 1$ and $\log(1+x) \ge x/(1+x)$ for $x \ge 0$, this is lower-bounded by $\rho^2\{1 - c(\delta) \cos^2 \tau\} \ge \rho^2 \sin^2 \tau$. Because $\sin(\pi - x) = -\sin x$,

$$(7.12) \qquad \inf_{\mathcal{N}_\gamma(\varepsilon) \times \Theta_\delta} \sin^2 \tau \ge \min\left\{\sin^2\left(\frac{\varepsilon}{2}\right), \sin^2(2\eta_4)\right\} > 0.$$

Since $\rho \ne 0$, (7.7) is proved.

Now consider (7.4). We can write $U_\beta(\theta) = \log Q(b_n(\beta))$, where $Q(s) = 1 + \hat{a}_1 s + \hat{a}_2 s^2$, $\hat{a}_1 = (\hat{g}_{11}^{(2)} \hat{g}_{22}^{(1)} - 2\hat{g}_{12}^{(1)} \hat{g}_{12}^{(2)})/\det\{G^{(1)}(\alpha)\}$, $\hat{a}_2 = (\hat{g}_{11}^{(3)} \hat{g}_{22}^{(1)} - \hat{g}_{12}^{(2)^2})/\det\{G_1^{(1)}(\alpha)\}$. For all $\theta$, $\hat{a}_2 \ge 0$ by the Cauchy inequality, and, since $\hat{\Omega}^*(\theta)$ and $\hat{G}^{(1)}(\alpha)$ are nonnegative definite, $Q(s)$ is nonnegative for all real $s$. It has a global minimum at $s = -\hat{a}_1/2\hat{a}_2$. Thus

$$\inf_{|s| \geq 1/\varepsilon} Q(s) \geq \left(1 - \frac{\hat{a}_1^2}{4\hat{a}_2}\right) 1\left(\left|\frac{\hat{a}_1}{2\hat{a}_2}\right| > \frac{1}{\varepsilon}\right) + \left(1 - \frac{|\hat{a}_1|}{\varepsilon} + \frac{\hat{a}_2}{\varepsilon^2}\right) 1\left(\left|\frac{\hat{a}_1}{2\hat{a}_2}\right| \leq \frac{1}{\varepsilon}\right)$$

$$= 1 - \frac{|\hat{a}_1|}{\varepsilon} + \frac{\hat{a}_2}{\varepsilon^2} + \left(\frac{|\hat{a}_1|}{\varepsilon} - \frac{\hat{a}_2}{\varepsilon^2} - \frac{\hat{a}_1^2}{4\hat{a}_2}\right) 1\left(\left|\frac{\hat{a}_1}{2\hat{a}_2}\right| > \frac{1}{\varepsilon}\right),$$

where $1(\cdot)$ denotes the indicator function. Thus the probability on the left-hand side of (7.4) is bounded by

$$
(7.13) \qquad
\begin{aligned}
& P\left(\log\left\{1 - \sup_{\Theta_\alpha} \frac{|\hat{a}_1|}{\varepsilon} + \inf_{\Theta_\alpha} \frac{\hat{a}_2}{\varepsilon^2}\right\} \leq 0\right) + P\left(\sup_{\Theta_\alpha} \left|\frac{\hat{a}_1}{2\hat{a}_2}\right| > \frac{1}{\varepsilon}\right) \\
& \leq 2P\left(\sup_{\Theta_\alpha} |\hat{a}_1 - a_1| + \frac{2}{\varepsilon} \sup_{\Theta_\alpha} |\hat{a}_2 - a_2| \geq \frac{1}{\varepsilon} \inf_{\Theta_\alpha} a_2 - \sup_{\Theta_\alpha} |a_1|\right)
\end{aligned}
$$

by elementary inequalities, where $a_1 = (g_{11}^{(2)} g_{22}^{(1)} - 2g_{12}^{(1)} g_{12}^{(2)})/\det G^{(i)}(\alpha)$, $a_2 = (g_{11}^{(3)} g_{22}^{(1)} - g_{12}^{(2)2})/\det G^{(i)}(\alpha)$. Now $\sup_{\Theta_\alpha} |\hat{g}_{k\ell}^{(i)} - g_{k\ell}^{(i)}| \to_p 0$ $(i = 2, 3,\ k, \ell = 1, 2)$ as $n \to \infty$, by the same method of proof as described for (7.5), so $\sup_{\Theta_\alpha} |\hat{a}_i - a_i| \to_p 0$ $(i = 1, 2)$ as $n \to \infty$. We need to show that the right-hand side of the last inequality in (7.13) is positive. It is easily seen that $\sup_{\Theta_\alpha} |a_1| < \infty$, noting boundedness away from zero on $\Theta_\alpha$ of denominators in the $g_{k\ell}^{(i)}$. Since $\varepsilon$ can be arbitrarily small we require only that $\inf_{\Theta_\alpha} a_2 > 0$. This is true because, on $\Theta_\alpha$,

$$
\begin{aligned}
g_{11}^{(3)} g_{22}^{(1)} - g_{12}^{(2)2} &= \omega_{22}^2 \left\{\frac{1}{\{2(\delta_1 - \delta_{02}) + 1\}(2\zeta_2 + 1)} - \frac{\cos^2 \gamma_0}{(\delta_1 + \zeta_1 - \delta_{02} + 1)^2}\right\} \\
&\geq \omega_{22}^2 \left[\frac{1}{\{2(\delta_1 - \delta_{02}) + 1\}(2\zeta_2 + 1)} - \frac{1}{(\delta_1 - \delta_{02} + \zeta_2 + 1)^2}\right] \\
&> \frac{\omega_{22}^2 \nu^2}{8} \geq \frac{\omega_{22}^2 \eta_2^2}{8} > 0.
\end{aligned}
$$

This completes the proof that $\hat{\alpha} \xrightarrow{p} \alpha_0$, $\hat{\beta} = \beta_0 + O_p((m/n)^{\nu_0})$. To replace "$O$" by "$o$" in the latter, for $\varepsilon \in (0, 1)$ define $\mathcal{N}^*(\varepsilon) = \mathcal{N}_\beta(\varepsilon^{1/2}(m/n)^{\nu_0}) \times \mathcal{N}_\gamma(\varepsilon^2) \times \mathcal{N}_\delta(\varepsilon^2)$, $\bar{\mathcal{N}}^*(\varepsilon) = \Theta \setminus \mathcal{N}^*(\varepsilon)$. We have $P(\hat{\theta} \in \bar{\mathcal{N}}^*(\varepsilon)) \leq P(\hat{\theta} \in \bar{\mathcal{N}}^*(\varepsilon) \cap \mathcal{N}(\varepsilon)) + P(\hat{\theta} \in \bar{\mathcal{N}}(\varepsilon))$. We have just shown that the last probability tends to zero. For the previous one it suffices to show that as $n \to \infty$

$$
(7.14) \qquad P\left(\inf_{\bar{\mathcal{N}}_\alpha^*(\varepsilon)} U_\alpha(\alpha) \leq 0\right) \to 0,
$$

$$
(7.15) \qquad P\left(\inf_{\bar{\mathcal{N}}_\beta(\varepsilon^{1/2}(m/n)^{\nu_0}) \times \mathcal{N}_\alpha(\varepsilon)} U_\beta(\theta) \leq 0\right) \to 0,
$$

where $\mathcal{N}_\alpha(\varepsilon) = \mathcal{N}_\gamma(\varepsilon) \times N_\delta(\varepsilon)$, $\bar{\mathcal{N}}_\alpha^*(\varepsilon) = \mathcal{N}_\alpha(\varepsilon) \setminus \mathcal{N}_\alpha^*(\varepsilon)$. The proof of (7.14) is as above. To prove (7.15), following the argument up to (7.13) we have to show

$$P\left( \sup_{\mathcal{N}_\alpha(\varepsilon)} |\hat{a}_1 - a_1| + 2\varepsilon^{1/2} \sup_{\mathcal{N}_\alpha(\varepsilon)} |\hat{a}_2 - a_2| \geq \varepsilon^{1/2} \inf_{\mathcal{N}_\alpha(\varepsilon)} a_2 - \sup_{\mathcal{N}_\alpha(\varepsilon)} |a_1| \right)$$
(7.16)
$$\to 0.$$

In view of the above remarks about (7.13) it remains to show that the right-hand side of the inequality in (7.16) is positive. We have

$$\sup_{\mathcal{N}_\alpha(\varepsilon)} |a_1| \leq \left\{ \sup_{\mathcal{N}_\alpha(\varepsilon)} \left| g_{11}^{(2)} g_{22}^{(1)} - 2 g_{12}^{(1)} g_{12}^{(2)} \right| \right\} \Big/ \inf_{\mathcal{N}_\alpha(\varepsilon)} \det\{G_1^{(1)}(\alpha)\}.$$

The denominator is already known to be finite and the quantity on the right-hand side whose absolute value is taken equals

$$2\omega_{12}\omega_{22}\left[ \frac{\cos\gamma_0}{(2\zeta_2 + 1)(\zeta_1 + \delta_1 - \delta_{02} - 1)} - \frac{\cos\gamma\cos\tau}{(\zeta_1 + \zeta_2 + 1)(\delta_1 + \zeta_2 - \delta_{02} + 1)} \right].$$

After rearrangement and application of trigonometric addition formula, this is seen to be bounded in absolute value by $K(|\gamma - \gamma_0| + \|\delta - \delta_0\|)$. It follows that $\sup_{\mathcal{N}_\alpha(\varepsilon)} |a_1| \leq K\varepsilon$. From the proof of Theorem 3, $\varepsilon^{1/2} \inf_{\mathcal{N}_\alpha(\varepsilon)} a_2 - \sup |a_1| \geq \varepsilon^{1/2}/K - K\varepsilon$, which, for arbitrarily large $K$, is bounded below by $\varepsilon^{1/2}/2K > 0$, choosing $\varepsilon \in (0, (4K^4)^{-1})$.

## 8. Proof of Theorem 4.
Define $s(\theta) = (\partial/\partial\theta)R(\theta)$, $S(\theta) = (\partial/\partial\theta')s(\theta)$. Denote by $\tilde{S}$ the matrix $S(\theta)$ when its $k$th row is evaluated at $\theta = \tilde{\theta}^{(k)}$. If $\|\tilde{\theta}^{(k)} - \theta_0\| \leq \|\hat{\theta} - \theta_0\|$, $k = 1, \ldots, 4$, the mean value theorem gives $\hat{\theta} - \theta_0 = \tilde{S}^{-1}s(\theta_0)$, for some such $\tilde{\theta}^{(k)}$. The theorem is established if

$$m^{1/2}\Delta_n^{-1}s(\theta_0) \xrightarrow{d} N_4(0, \Sigma), \tag{8.1}$$

$$\Delta_n^{-1}\tilde{S}\Delta_n^{-1} \xrightarrow{p} \Sigma. \tag{8.2}$$

Denoting by $\theta_k$, $s_k(\theta)$, the $k$th elements of $\theta$, $s(\theta)$, and by $s_{k\ell}(\theta)$ the $(k, \ell)$th element of $S(\theta)$,

$$s_k(\theta) = \text{tr}\left\{ \frac{\partial\hat{\Omega}(\theta)}{\partial\theta_k} \hat{\Omega}(\theta)^{-1} \right\} - 1(k = 3 \text{ or } 4) \frac{2}{m} \sum_j \log|\psi(\lambda_j; \gamma)|, \tag{8.3}$$

$$s_{k\ell}(\theta) = \text{tr}\left\{ \frac{\partial^2\hat{\Omega}(\theta)}{\partial\theta_k \partial\theta_\ell} \hat{\Omega}(\theta)^{-1} - \frac{\partial\hat{\Omega}(\theta)}{\partial\theta_k} \hat{\Omega}(\theta)^{-1} \frac{\partial\hat{\Omega}(\theta)}{\partial\theta_\ell} \hat{\Omega}(\theta)^{-1} \right\}. \tag{8.4}$$

Now

$$\frac{\partial\hat{\Omega}(\theta)}{\partial\theta_k} = \text{Re}\left\{ \frac{1}{m} \sum_j A_j^{(k)} \right\}, \qquad \frac{\partial^2\hat{\Omega}(\theta)}{\partial\theta_k \partial\theta_\ell} = \text{Re}\left\{ \frac{1}{m} \sum_j A_j^{(k,\ell)} \right\},$$

$k, \ell = 1, 2$, writing $A_j^{(k)} = (\partial/\partial\theta_k)A_j$, $A^{(k,\ell)} = (\partial/\partial\theta_\ell)A_j^{(k)}$, $A_j = A(\lambda_j;\theta)$. To simplify we proceed, as in the proof of Theorem 3, as if $\psi(\lambda) = |\lambda|$. This can be justified via Assumption B3; further discussion appears later in the proof. Define $E_{k\ell}$ by replacing the $(k,\ell)$th element by 1 in the $2\times 2$ matrix of zeros. Noting that $E_{12}B' = -E_{12}$ we deduce $A_j^{(1)} = -\lambda_j^{-\nu}(E_{12}A_j e^{i\gamma} - A_j E_{21}e^{-i\gamma})$, $A_j^{(2)} = iA_jE_{22} - iE_{22}A_j$, $A_j^{(2+k)} = (\log\lambda_j)(E_{kk}A_j + A_jE_{kk})$, $k = 1, 2$, $A_j^{(1,1)} = 2\lambda_j^{-2\nu}E_{12}A_jE_{21}$, $A_j^{(1,2)} = i\lambda_j^{-\nu}(E_{22}A_jE_{21} - E_{12}A_jE_{22})$, $A_j^{(1,2+k)} = -(\log\lambda_j)\times \lambda_j^{-\nu}(E_{kk}E_{12}A_j + E_{kk}A_jE_{21} + E_{12}A_jE_{kk} + A_jE_{21}E_{kk})$, $A_j^{(2,2)} = 2E_{22}A_jE_{22} - E_{22}A_j - AE_{22}$, $A_j^{(2,2+k)} = -i(\log\lambda_j)(E_{kk}A_jE_{22} - E_{kk}E_{22}A_j - E_{22}A_jE_{kk} + A_jE_{22}E_{kk})$, $A_j^{(2+k,2+\ell)} = (\log\lambda_j)^2(E_{kk}E_{\ell\ell}A_j + A_j \times E_{\ell\ell}E_{kk} + E_{kk}A_jE_{\ell\ell} + E_{\ell\ell}A_jE_{kk})$. Thus from (8.3), with $A_{0j} = A(\lambda_j;\theta_0)$,

$$s_1(\theta_0) = -\operatorname{tr}\frac{1}{m}\sum_j \lambda_j^{-\nu_0}(E_{12}A_{0j}e^{i\gamma_0} + A_{0j}E_{21}e^{-i\gamma_0})\hat\Omega(\theta_0)^{-1},$$

$$s_2(\theta_0) = i\operatorname{tr}\left\{\frac{1}{m}\sum_j(A_{0j}E_{22} - E_{22}A_{0j})\hat\Omega(\theta_0)^{-1}\right\},$$

$$s_{2+k}(\theta_0) = \operatorname{tr}\frac{1}{m}\sum_j\left(\log\lambda_j - \frac{1}{m}\sum_i\log\lambda_i\right)(E_{kk}A_{0j} + A_{0j}E_{kk})\hat\Omega(\theta_0)^{-1},$$

for $k = 1, 2$, where the real part operator is omitted because imaginary parts are automatically eliminated here, and we use $\hat\Omega(\theta_0) = m^{-1}\sum_j\operatorname{Re}\{A_{0j}\}$. We can replace, with negligible error, $\hat\Omega(\theta_0)$ by $\Omega_0$ and $A_{0j}$ by $T_j = PI_{\varepsilon j}P'$ in $m^{1/2}\Delta_n^{-1}s(\theta_0)$, using arguments of [22, 29, 33], and allowing $p \geq 2$. Thus $m^{1/2}\Delta_n^{-1}s(\theta_0)$ differs by $o_p(1)$ from $m^{1/2}\Delta_n^{-1}s^*(\theta_0)$, where $s^*(\theta_0)$ has $k$th element

$$(8.5)\qquad s_k^* = \frac{2}{m}\sum_j\operatorname{tr}(U_{Rkj}\operatorname{Re}\{I_{\varepsilon j}\} + U_{Ikj}\operatorname{Im}\{I_{\varepsilon j}\}),$$

where $U_{R1j} = -\cos\gamma_0(\lambda_j^{-\nu_0} - m^{-1}\sum_i\lambda_i^{-\nu_0})P'\Omega_0^{-1}E_{12}P$, $U_{I1j} = -\sin\gamma_0\lambda_j^{-\nu_0}P'\Omega_0^{-1}E_{12}P$, $U_{R2j} = 0$, $U_{I2j} = P'\Omega_0^{-1}E_{22}P$, $U_{R,2+k,j} = (\log j - m^{-1}\sum_i\log i)\times P'\Omega_0^{-1}E_{kk}P$, $U_{I,2+k,j} = 0$, for $k = 1, 2$. After rearrangement and application of a martingale CLT we deduce, following the same references, $\Delta_n^{-1}s^*/2 \to_d N_4(0, \Sigma)$. [The formula for $\Sigma$ can be most easily verified after noting that $Es_k^*s_\ell^* = 8m^{-2}\sum_j\operatorname{tr}\{U_{Rkj}(U_{R\ell j}' + U_{R\ell j}) + U_{Ikj}(U_{I\ell j}' - U_{I\ell j})\}$, plus a negligible fourth cumulant term.] This completes the proof of (8.1).

Turning to (8.2), it suffices to show that

$$(8.6)\qquad \Delta_n^{-1}\{\tilde S - S(\theta_0)\}\Delta_n^{-1} \xrightarrow{p} 0,$$

$$(8.7) \qquad \tfrac{1}{2}\Delta_n^{-1} S(\theta_0) \Delta_n^{-1} \xrightarrow{p} \Sigma.$$

We omit the straightforward proof of (8.7). To prove (8.6), we require a rate of convergence for the $\tilde{\delta}_i$. Put $\tilde{\theta} = (\bar{\beta}, \tilde{\alpha}')' = (\bar{\beta}, \tilde{\gamma}, \tilde{\delta}_1, \tilde{\delta}_2)'$, for $\|\tilde{\theta} - \theta_0\| \leq \|\hat{\theta} - \theta_0\|$. For some such $\tilde{\theta}$, $\hat{\Omega}(\tilde{\theta})$ appears in all elements of $\tilde{S}$. From Section 7 we can write, with the same definitions, and $\Psi$ again replaced by $\Phi$, $\hat{\Omega}(\tilde{\theta}) = \Xi(\tilde{\delta})\{\hat{G}^{(1)}(\tilde{\alpha}) + b_n(\tilde{\beta})\hat{G}^{(2)}(\tilde{\alpha}) + b_n^2(\tilde{\beta})\hat{G}^{(3)}(\tilde{\alpha})\}\Xi(\tilde{\delta})$. Then from Theorem 3, $\hat{\Omega}(\tilde{\theta}) - \hat{\Omega}(\theta_0) \to_p 0$ if $\tilde{\delta} \to_p \delta_0$ and $\hat{G}^{(i)}(\tilde{\alpha}) - \hat{G}^{(i)}(\alpha_0) \to_p 0$, $i = 1, 2, 3$. To achieve the latter, $H_j = A_{0j}$ can be replaced as before by the $T_j$, but from the definitions of Section 7 the $\tilde{\delta}_k$ are involved as exponents of $(j/n)$, $j = 1, \ldots, m$, in the $\hat{G}^{(1)}(\tilde{\alpha})$, so more than the consistency established in Theorem 3 is needed (though consistency of $\tilde{\gamma}$ suffices). So far as remaining terms which make up elements of $\tilde{S}$ are concerned, similar considerations apply, indeed differentiation produces factors $\log|\psi_j|$, $\log^2|\psi_j|$ in some summands. In [29], only $\psi(\lambda) = |\lambda|$ was considered, and $\log n$ terms are precisely eliminated prior to taking limits, as in Section 7. With more general $\psi_j$ this does not happen, as in [33]'s choice of $\psi$, and as there we establish something a little stronger. It suffices to show that $(\log n)^C (\hat{\delta}_k - \delta_{0k}) \to_p 0$, $k = 1, 2$, for any $C < \infty$ [explaining the requirement $(\log n)^C/m \to 0$ in (3.7)]. Arguing as before, this follows if, as $n \to \infty$, $\sup_{\mathcal{N}_\alpha(\varepsilon)} \|\Upsilon(\delta)\{\hat{G}^{(1)}(\alpha) - G^{(1)}(\alpha)\}\Upsilon(\delta)\| = o_p((\log n)^{-2C})$, $\inf_{\bar{\mathcal{N}}_{\tilde{\gamma}}(\varepsilon^2) \times \mathcal{N}_\delta(\varepsilon)} \log \det\{\Upsilon(\delta)G^{(1)}(\alpha)\Upsilon(\delta)G^{(1)}(\alpha_0)^{-1}\} > 0$, $\varliminf_{n \to \infty}(\log n)^{2C} \times \inf_{\bar{\mathcal{N}}_\delta(\varepsilon/(\log n)^C)} u(\delta) > 0$, for any $\varepsilon \in (0, 1)$. The first result follows by straightforward extension of the proof of (4.6) in [29], the rate being due to $\varepsilon_t$ now having a finite moment of order greater than 2. The proof of the second is identical to that of (7.7), the only difference in outcome being the replacement of $\varepsilon$ by $\varepsilon^2$. As in the proof of [29], we deduce the final result from the inequalities in (7.9).

**Acknowledgments.** I thank the Associate Editor and three referees for comments which have led to significant improvements, and Afonso Goncalves da Silva for help with the numerical work.

## REFERENCES

[1] ADENSTEDT, R. (1974). On large-sample estimation for the mean of a stationary random sequence. *Ann. Statist.* **2** 1095–1107. MR0368354

[2] ANDERSON, T. W. and WALKER, A. M. (1964). On the asymptotic distribution of the autocorrelations of a sample from a linear stochastic process. *Ann. Math. Statist.* **35** 1296–1303. MR0165602

[3] BREIDT, F. J., DAVIS, R. A. and TRINDADE, A. A. (2001). Least absolute deviation estimation for all-pass time series models. *Ann. Statist.* **29** 919–946. MR1869234

[4] BRILLINGER, D. R. (1975). *Time Series. Data Analysis and Theory.* Holt, Rinehart and Winston, New York. MR0443257

[5] CHEN, W. W. and HURVICH, C. M. (2006). Semiparametric estimation of fractional cointegrating subspaces. *Ann. Statist.* **34** 2939–2979. MR2329474

[6] CHRISTENSEN, B. J. and NIELSEN, M. (2006). Asymptotic normality of narrow-band least squares in the stationary fractional cointegration model and volatility forecasting. *J. Econometrics* **133** 343–371. MR2250183

[7] DAHLHAUS, R. (1989). Efficient parameter estimation for self-similar processes. *Ann. Statist.* **17** 1749–1766. MR1026311

[8] DAVIDSON, J. E. H. and HASHIMADZE, N. (2008). Alternative frequency and time domain versions of fractional Brownian motion. *Econometric Theory* **24** 256–297.

[9] FOX, R. and TAQQU, M. S. (1986). Large-sample properties of parameter estimates for strongly dependent stationary Gaussian series. *Ann. Statist.* **14** 517–532. MR0840512

[10] GIRAITIS, L. and ROBINSON, P. M. (2003). Edgeworth expansions for semiparametric Whittle estimation of long memory. *Ann. Statist.* **31** 1325–1375. MR2001652

[11] GIRAITIS, L. and SURGAILIS, D. (1990). A central limit theorem for quadratic forms in strongly dependent linear variables and its application to asymptotic normality of Whittle's estimate. *Probab. Theory Related Fields* **86** 87–104. MR1061950

[12] HANNAN, E. J. (1970). *Multiple Time Series.* Wiley, New York. MR0279952

[13] HENRY, M. and ROBINSON, P. M. (1996). Bandwidth choice in Gaussian semiparametric estimation of long range dependence. *Athens Conference on Applied Probability and Time Series Analysis* **II**. *Lecture Notes in Statist.* **115** 220–232. Springer, New York. MR1466748

[14] HOSOYA, Y. (1997). Limit theory with long-range dependence and statistical inference of related models. *Ann. Statist.* **25** 105–137. MR1429919

[15] HUALDE, J. and ROBINSON, P. M. (2004). Semiparametric estimation of fractional cointegration. Preprint, London School of Economics.

[16] HURVICH, C. M. and BELTRAO, K. J. (1994). Automatic semiparametric estimation of the memory parameter of a long memory time series. *J. Time Ser. Anal.* **15** 285–302. MR1278413

[17] HURVICH, C. M. and DEO, R. S. (2001). Plug-in selection of the number of frequencies in regression estimates of the memory parameter of long memory time series. *J. Time Ser. Anal.* **20** 331–341. MR1693157

[18] HURVICH, C. M., MOULINES, E. and SOULIER, E. (2005). Estimating long memory in volatility. *Econometrica* **73** 1283–1328. MR2149248

[19] JOHANSEN, S. (2005). A representation theory for a class of vector autoregressive models for fractional processes. Working paper, Dept. Applied Mathematics and Statistics, Univ. Copenhagen.

[20] KÜNSCH, H. R. (1987). Statistical aspects of self-similar processes. In *Proceedings of the First World Congress of the Bernoulli Society* **1** 67–74. VNU Sci. Press, Utrecht. MR1092336

[21] LII, K. S. and ROSENBLATT, M. (1982). Deconvolution and estimation of transfer function phase and coefficients for non-Gaussian linear processes. *Ann. Statist.* **10** 1195–1208. MR0673654

[22] LOBATO, I. N. (1999). A semiparametric two-step estimator in a multivariate long memory model. *J. Econometrics* **90** 129–153. MR1682397

[23] NIELSEN, M. O. (2005). Semiparametric estimation in time-series regression with long-range dependence. *J. Time Ser. Anal.* **26** 279–304. MR2122898

[24] NIELSEN, M. O. (2007). Local Whittle analysis of stationary fractional cointegration and the implied-realized volatility relation. *J. Bus. Econ. Statist.* **25** 427–446.

[25] PARZEN, E. (1957). A central limit theorem for multilinear stochastic process. *Ann. Math. Statist.* **28** 252–256. MR0084899

[26] PHILLIPS, P. C. B. and SHIMOTSU, K. (2004). Whittle estimation in nonstationary and unit root cases. *Ann. Statist.* **32** 656–672. MR2060173

[27] ROBINSON, P. M. (1994). Semiparametric analysis of long-memory time series. *Ann. Statist.* **22** 515–539. MR1272097

[28] ROBINSON, P. M. (1995). Log-periodogram regression of time series with long range dependence. *Ann. Statist.* **23** 1048–1072. MR1345214

[29] ROBINSON, P. M. (1995). Gaussian semiparametric estimation of long range dependence. *Ann. Statist.* **23** 1630–1661. MR1370301

[30] ROBINSON, P. M. (2005). Robust covariance matrix estimation: "HAC" estimates with long memory/antipersistence correction. *Econometric Theory* **21** 171–180. MR2153861

[31] ROBINSON, P. M. and MARINUCCI, D. (2001). Narrow-band analysis of nonstationary processes. *Ann. Statist.* **29** 947–986. MR1869235

[32] ROBINSON, P. M. and MARINUCCI, D. (2003). Semiparametric frequency-domain analysis of fractional cointegration. In *Time Series with Long Memory. Adv. Texts Econometrics* 334–373. Oxford Univ. Press. MR2077521

[33] SHIMOTSU, K. (2007). Gaussian semiparametric estimation of multivariate fractionally integrated processes. *J. Econometrics* **137** 277–310. MR2354946

[34] SHIMOTSU, K. and PHILLIPS, P. C. B. (2006). Exact local Whittle estimation of fractional integration. *Ann. Statist.* **33** 1890–1933. MR2166565

[35] VELASCO, C. (1999). Gaussian semiparametric estimation of non-stationary time series. *J. Time Ser. Anal.* **20** 87–127. MR1678573

[36] VELASCO, C. (2000). Gaussian semi-parametric estimation of fractional cointegration. Preprint, Univ. Carlos III de Madrid.

[37] VELASCO, C. (2003). Gaussian semi-parametric estimation of fractional cointegration. *J. Time Ser. Anal.* **24** 345–378. MR1984601

[38] ZYGMUND, A. A. (1977). *Trigonometric Series.* **I**, **II**. Cambridge Univ. Press, New York. MR0617944

DEPARTMENT OF ECONOMICS
LONDON SCHOOL OF ECONOMICS
HOUGHTON STREET
LONDON WC2A 2AE
UNITED KINGDOM
E-MAIL: p.m.robinson@lse.ac.uk